\documentclass{article}
\usepackage{amsfonts}
\title{\large \bf Warped product contact $CR$-submanifolds of globally framed $f$-manifolds with Lorentz metric }
\author{\small \bf Khushwant Singh and S. S. Bhatia \\
\small \em School of Mathematics and Computer Applications\\
\small \em Thapar University, Patiala - 147 004,
India\\
\small E-mail: khushwantchahil@gmail.com, ssbhatia@thapa.redu}

\date{}
\newtheorem{thm}{Theorem}[section]
\newtheorem{lm}[thm]{Lemma}
\newtheorem{pp}[thm]{Proposition}

\newtheorem{cor}[thm]{Corollary}

\begin{document}

\maketitle
\begin{abstract}
In the present paper, we study globally framed $f$-manifolds in the
particular setting of indefinite $S$-manifolds for both spacelike
and timelike cases. We prove that if $M = N^\perp\times_f N^T$ is a
warped $CR$-submanifold such that $N^\perp$ is $\phi$-anti-invariant
and $N^T$ is $\phi$-invariant, then $M$ is a $CR$-product. We show
that the second fundamental form of a contact $CR$ warped product of
a indefinite $S$ space form satisfies a geometric inequality, $
\|h\|^2\geq p\{3\|\nabla ln f\|^2-\triangle ln f+(c+2)k+1 \}.$
\end{abstract}

\vspace{2ex}

\noindent {\bf 2000 Mathematics Subject Classification:}
53C25,53C40, 53B25.

\vspace{2ex}

\noindent {\bf Keywords:} Indefinite $S$-manifold, Globally framed
$f$-manifolds, Warped product, Contact $CR$-submanifolds.

\section{\large Introduction}

\par It is well-known that the notion of warped products plays some important role in differential
geometry as well as physics. R.L. Bishop and B. O'Neill in 1969
introduced the concept of a warped product manifold to provide a
class of complete Riemannian manifolds with everywhere negative
curvature. The warped product scheme was later applied to
semi-Riemannian geometry ([1],[2]) and general relativity [3].
Recently, Chen [10] (see also [11]) studied warped product
$CR$-submanifolds and showed that there exist no warped product
$CR$-submanifolds of the form $M=N^\perp\times{_{f}N^T}$ such that
$N^\perp$ is a totally real submanifold and $N^T$ is a holomorphic
submanifold of a Kaehler manifold $\tilde {M}$. Therefore he
considered warped product $CR$-submanifold in the form
$M=N^T\times{_{f}N^\perp}$ which is called $CR$-warped product,
where $N^T$ and $N^\perp$ are holomorphic and totally real
submanifolds of a Kaehler manifold $\tilde{M}$. Motivated by Chen's
papers many authors studied
$CR$-warped product submanifolds in almost complex as well as contact setting (see [8], [13]).\\

\section{Globally framed $f$-manifolds with Lorentz metric}

\noindent First we recall some definitions due ([1],[6],[7],[9]).\\

\noindent A $(2n+1)-$dimensional $C^\infty$ manifold $\bar M$ is
said to have an {\it{almost contact structure}} if there exist on
$\bar M$ a tensor field $\phi$ of type $(1,1)$, a vector field $\xi$
and a $1-$form $\eta$ satisfying
$$\phi^2=-I+\eta\otimes\xi,~~\phi\xi=0,~~\eta\circ\phi=0,~~\eta(\xi) = 1.$$
There always exists a Riemannian metric $g$ on an almost contact
manifold $\bar M$ satisfying the following compatibility condition
$$\eta(X)=g(X, \xi),~~~g(\phi X, \phi Y)=g(X, Y)-\eta(X)\eta(Y)$$
where $X$ and $Y$ are vector fields on $\bar M$.\\

\noindent Moreover, if $g$ is a semi-riemannian metric on
$\tilde{M}^{2n+1}$ such that,
$$g(\phi X,\phi Y)=g(X,Y)-\epsilon \eta(X)\eta(Y), $$

\noindent where $\epsilon=\pm1$ according to the causal character of
$\xi$, $\tilde{M}^{2n+1}$ is called an indefinite almost contact
manifold if $d\eta=\Phi$, $\Phi$ being defined by $\Phi
(X,Y)=g(X,fY)$.\\

\noindent In the Riemannian case a generalization of these
structures have  been studied by Blair [6], by Goldberg and Yano
[12]. Brunetti and Pastore [9] studied such structures in
semi-Riemannian context.\\

\noindent  A manifold $\tilde{M}$ is called a globally framed
$f$-manifold (briefly $g.f.f$-manifold) if it is endowed with a non
null $(1,1)-$tensor field $\phi$ of constant rank, such that $ker~
\phi$ is parallelizable i.e. there exist global vector field
$\xi_\alpha$, such that $\alpha=\{1,...,s\}$, and $1$-form
$\eta^{\alpha}$, satisfying
$$\phi^2=-I+\eta^{\alpha}\otimes\xi_{\alpha}~~ and~~ \eta^{\alpha}(\xi_\beta)={\delta^{\alpha}}_\beta. \leqno(2.1) $$

\noindent A $g.f.f$-manifold
$(\tilde{M}^{2n+s},\phi,\eta^\alpha,\xi_\alpha)$, such that
$\alpha\in\{1,...,s\}$, is said to be an indefinite $g.f.f$-manifold
if $g$ is a semi-Riemannian metric satisfying the following
compatibility condition

$$ g(\phi X,\phi Y)=g(X,Y)-\epsilon_{\alpha} \eta^{\alpha}(X)\eta^{\alpha} (Y), \leqno(2.2) $$

\noindent for any vector fields $X,Y$, where $\epsilon_\alpha=\pm1$
according to whether $\xi_\alpha$ is spacelike or timelike. Then for
any $\alpha\in \{1,...,s\}$,
$$\eta^{\alpha} (X)=\epsilon_\alpha g(X,\xi_\alpha).$$

\noindent {\bf Note:} We will consider $\alpha\in \{1,...,s\}$
throughout the paper.\\

\noindent An indefinite $g.f.f$-manifold is an indefinite
$S$-manifold if it is normal and $d\eta^\alpha=\Phi$, for any
$\alpha\in\{1,...,s\}$, where $\Phi(X,Y)=g(X,\phi Y)$ for any $X,Y
\in \chi(M^{2n+s})$. The Levi-Civita connection of an indefinite
$S$-manifold satisfies

$$(\nabla_X \phi)Y=g(\phi X,\phi Y)\bar \xi_\alpha+\bar\eta^\alpha(Y)\phi^2 (X), \leqno(2.3) $$

\noindent where $\tilde{\xi}=\sum\limits_{\alpha=1}\limits^s
\xi_\alpha$ and $\tilde{\eta}=\epsilon_\alpha \eta^\alpha$. Note
that for $s=1$, indefinite $S$-manifold becomes indefinite Sasaki manifold.\\

\noindent From (2.3), it follows that $\nabla_X
\xi_\alpha=-\epsilon_\alpha \phi X$ and $ker \phi$ is an integrable
flat distribution since $\nabla_{\xi_\alpha} \xi_\beta=0$, for any
$\alpha,\beta\in\{1,...,s\}$. A $g.f.f$-manifold is subject to the
topological condition: It has to be either non compact or compact
with vanishing Euler characteristic, since it admits never vanishing
vector fields. This implies that such a $g.f.f$-manifold always admit Lorentz metrics.\\

\noindent An indefinite $S$-manifold
$(\tilde{M},\phi,\xi_\alpha,\eta^\alpha,g)$ is said to be an
indefinite $S-$space if the $\phi$-sectional curvature $H_p(X)$ is
constant, for any point and any $\phi$-plane. In particular, in [9]
it is proved that an indefinite $S$-manifold
$(M,\phi,\xi_\alpha,\eta^\alpha,g)$ is an indefinite $S-$space form
with $H_p (X)=c$ if and only if the Riemannian $(0,4)-$type
curvature tensor field $R$ is given by
$$R(X,Y,Z,W)=-\frac{c+3\epsilon}{4}\{g(\phi Y,phi Z)g(\phi X,\phi W)-g(\phi X,\phi Z)g(\phi Y,\phi W)\} \leqno(2.4) $$
$$~~~~~~~-\frac{c-\epsilon}{4} \{\Phi(W,X)\Phi(Z,Y)-\Phi(Z,X)\Phi(W,Y) $$
$$~~~~~~~~~+2\Phi(X,Y)\Phi(W,Z) \}-\{\tilde{\eta}(W)\tilde{\eta}(X)g(\phi Z,\phi Y)$$
$$~~~~~~~~~~-\tilde{\eta}(W)\tilde{\eta}(Y)g(\phi Z,\phi X)+\tilde{\eta}(Y)\tilde{\eta}(Z)g(\phi W,\phi X)$$
$$-\tilde{\eta}(Z)\tilde{\eta}(X)g(\phi W,\phi Y)\},~~~~~~~~~~~~~~~~~ $$
\noindent for any vector fields $X$,$Y$,$Z$ and $W$ on $M$, where
$\epsilon=\sum\limits_{\alpha=1}\limits^s \epsilon_\alpha$.\\

\noindent Let $M$ be a real $m-$dimensional submanifold of $\tilde
{M}^{2n+s}$, tangent to the global vector field $\xi_\alpha$. We
shall need the Gauss and Weigarten formulae
$$\tilde {\nabla}_X Y=\nabla_X Y+ h(X,Y),~~\tilde {\nabla}_X N=-A_N X+ \nabla_X^\perp N, \leqno(2.5)$$
for any $X,Y\in \chi(M)$ and $N\in \Gamma^\infty(T(M)^\perp)$, where
$\nabla^\perp$ is the connection on the normal bundle $T(M)^\perp $,
$h$ is the second fundamental form and $A_N$ is the Weingarten map
associated with the vector field $N\in T(M)^\perp$ as
$$g(A_N X,Y)=\tilde{g}(h(X,Y),N). $$

\noindent For any $X\in\chi(M)$ we set $PX=tan(\phi X)$ and
$FX=nor(\phi X)$, where $tan_x$ and $nor_x$ are the natural
projections associated to the direct sum decomposition
$$ T_x(\tilde {M})=T_x (M)\oplus T{(M)_x}^\perp,~~~x\in M. $$

\noindent Then $P$ is an endomorphism of the tangent bundle of
$T(M)$ and $F$ is a normal bundle valued $1-$form on $M$. Since
$\xi_\alpha$ is tangent to $M$, we get\\

$$ P\xi_\alpha=0,~~~~F\xi_\alpha=0,~~~~\nabla_X\xi_\alpha=PX,~~~~h(X,\xi_\alpha)=FX. $$

\noindent Similarly, for a normal vector field $F$, we put
$tF=tan(\phi F)$ and $fF=nor(\phi F)$ for the tangential and normal
part of $\phi F$, respectively.\\

\noindent The covariant derivative of the morphisms $P$ and $F$ are
defined respectively as
$$(\nabla_U P)V=\nabla_U PV- P\nabla_U V,~~~~(\nabla_U F)V=\nabla^\perp_U FV-F \nabla_U V,$$
for $U, V\in \chi(M)$. On using equation (2.3) and (2.5) we get
$$(\nabla_U P)V=g(PU,PV)\xi_\alpha+\eta^\alpha(U)\eta^\alpha(V)\xi_\alpha-\eta^\alpha(V)U-th(U,V)-A_{FV} U \leqno(2.6) $$
\noindent and
$$ (\nabla_U F)V=g(FU,FV)\xi_\alpha-fh(U,V)-h(U,PV). \leqno(2.7) $$

\noindent The Riemannian curvature tensor $R$ of $M$ is given by

$$~~~R_{XY} Z=\frac{3\epsilon^2+c\epsilon-4}{4}\{g(Y,Z)\eta^\alpha(X)\xi_\alpha-g(X,Z)\eta^\alpha(Y) \xi_\alpha  \leqno(2.8) $$
$$~~~~~~~~~~~~~~~~~~~~~~~~~~~~~~~~~~~~~~~~~~+\eta^\alpha(Y)\eta^\alpha(Z)X -\eta^\alpha(X)\eta^\alpha(Z)Y+\frac{c+3\epsilon}{4}\{g(X,Z)Y-g(Y,Z)X\}  $$
$$~~~~~~~~~~~~~~~~~~~~~~~~~~-\frac{c-\epsilon}{4}\{g(Z,PY)PX-g(Z,PX)PY-2g(X,PY)PZ\}, $$

\noindent for all $X$,$Y$,$Z$ vector fields on $M$, we recall the
equation of Gauss and Codazzi, respectively

$$~~~~~~~~~~~~~~~~~~~~~\tilde{g}(\tilde{R}_{XY} Z,W)=g(R_{XY} Z,W)-g(h(X,W),h(Y,Z))+g(h(Y,W),h(X,Z)) \leqno(2.9)  $$
$$ (\tilde{\nabla}_X h)(Y,Z)-(\tilde{\nabla}_Y h)(X,Z)=(\tilde{R}_{XY} Z)^\perp, \leqno(2.10) $$

\noindent where $(\nabla)h$ the covariant derivative of the second
fundamental form is given by

$$ (\tilde{(\nabla}_X h)(Y,Z)=\bar\nabla_X h(Y,Z)-h(\nabla_X Y,Z)-h(Y,\nabla_X Z), \leqno(2.11) $$

\noindent for all $X,Y,Z \in TM$. Codazzi equation becomes
$$(\tilde{\nabla}_X h)(Y,Z)-(\tilde{\nabla}_Y h)(X,Z)=\frac{c-\epsilon}{4}\{g(PX,Z)FY-g(PY,Z)FX+2g(X,PY)FZ\}. \leqno(2.12) $$

\noindent The second fundamental form $h$ satisfies the classical
Codazzi equation (according to [3]) if

$$(\tilde{\nabla}_X h)(Y,Z)=(\tilde{\nabla}_Y h)(X,Z). \leqno(2.13) $$

\begin{lm} Let $M^m$ be a submanifold of an indefinite $S$-space form
$\tilde{M}^{2m+s} (c)$ tangent to the global vector field
$\xi_\alpha$ with $c\neq1,-1$  according to whether $\xi_\alpha$ is
spacelike or timelike. If the second fundamental form $h$ of $M^m$
satisfies the classical Codazzi equation then $M^m$ is
$\phi-$invariant or $\phi-$anti-invariant.
\end{lm}

\noindent{\bf{Proof.}} By using (2.12) and (2.13), we get
$$\{g(PX,Z)FY-g(PY,Z)FX+2g(X,PY)FZ\}=0, \leqno(2.14) $$
\noindent for all $X,Y,Z\in T(M)$. By contradiction, let there exist
$U_x\in T_x(M)$ such that $PU_x\neq0$ and $FU_x\neq0$. From (2.14),
we deduce $2g(U_x,PU_x)FU_x=0$, which is false. Therefore, for
$U_x\in T_x(M)$, we have either $PU_x=0$ or $FU_x=0$. It can be also
proved that we can not have $U_x,V_x\in T_x(M)$ such that
$PU_x\neq0$, $FU_x=0$, $PV_x=0$ and $FV_x\neq0$. Therefore either
$P=0$ or $F=0$ which completes the statement.\\

\begin{lm} Let $M^m$ be a contact $CR$-submanifold of an indefinite $S$-manifold $\tilde {M}^{2n+s}$. Then for any $Z,W\in D^\perp$, we have
$$A_{FZ} W+A_{FW}Z=\eta^\alpha(Z)W+\eta^\alpha(W)Z-2\eta^\alpha(W)\eta^\alpha(Z)\xi_\alpha. \leqno(2.15)  $$
\end{lm}
\noindent{\bf Proof.} Proof is straightforward and can be obtained
by using equation (1.3) and (1.5).\\

\noindent Clearly, for $\xi_\alpha\in D$ we also have
$$A_{FZ} W+A_{FW}Z=0, \leqno(2.16) $$
\noindent where $W,Z\in D^\perp$.\\

\begin{lm} Let $M^m$ be a contact $CR$-submanifold of an indefinite
$S$-manifold $\tilde {M}^{2n+s}$ with $\xi_\alpha\in D$. Then the
following are equivalent
\begin{enumerate}
 \item[{(i)}] $h(X,PY)=h(PX,Y)~~~~\forall~ X,Y\in D$,
 \item[{(ii)}] $\tilde{g}(h(X,PY),\phi Z)=\tilde{g}(h(PX,Y),\phi Z)~~~~\forall~ X,Y\in D,~~~~\forall~ Z\in D^\perp$,
 \item[{(iii)}] $D$ is completely Integrable.
\end{enumerate}
\end{lm}
\noindent{\bf Proof.} Easy Calculations.\\

\noindent For a leaf of anti-invariant distribution $D^\perp$. We
prove the following

\begin{pp} Let $M^m$ be a contact $CR$-submanifold of an indefinite
$S$-manifold $\tilde M^{2n+s}$. Then any leaf of $D^\perp$ is
totally geodesic in $M^m$ if and only if
$$g(h(D,D^\perp),\phi D^\perp)=0. \leqno(2.17) $$
\end{pp}
\noindent {\bf{Proof.}} By hypothesis
$$g(P\nabla_W Z,Y)=-g(A_{FZ} W,Y)-g(th(Z,W),Y)=-g(h(Y,W),FZ), \leqno(2.18) $$
\noindent for any $Y\in D$, $Z,W\in D^\perp$. Then
$$g(\nabla_W Z,PY)=-g(h(Y,W),FZ). \leqno(2.19) $$

\noindent Consequently, $\nabla_W Z\in D^\perp$ if and
only if (2.17) holds.\\

\noindent Let $\nu$ be the complementary orthogonal subbundle of
$\phi D^\perp$ in the normal bundle $T(M)^\perp$. Thus we have the
following direct sum decomposition
$$ T(M)^\perp=\phi D^\perp \oplus \nu. \leqno(2.20) $$

\noindent Similarly, we can also prove that, $h(X,Y)\in \nu$ and
$\phi h(X,Y)=h(X,PY)$ for all $X,Y$ tangent to $N^T$. On $N^T$ we
have an induced indefinite $S$-structure.\\

\begin{lm} Let $M^m$ be a contact $CR$-submanifold of an indefinite
$S$-manifold $\tilde {M}^{2n+s}$ with $\xi_\alpha\in D$. Then for
all $X,Y\in D$, we have $\phi h(X,Y)\in \phi D^\perp\oplus\nu.$
\end{lm}
\noindent{\bf Proof.} From lemma 2.5 it follows that $\phi \nu=\nu$.
Since $h$ is normal to $M^m$ and $\eta^\alpha(D^\perp)=0$.
We easily get the result.\\

\noindent A contact $CR$-submanifold $M^m$ of an indefinite
$S$-manifold is called contact $CR$-product if it is locally a
Riemannian product of a $\phi-$invariant submanifold $N^T$
tangent to $\xi_\alpha$ and a totally real submanifold $N^\perp$ of $\tilde{M}^{2n+s}$.\\

\begin{thm} Let $M^m$ be a contact $CR$-submanifold of an indefinite
$S$-manifold $\tilde{M}^{2n+s}$ and set $\xi_\alpha \in D$. Then
$M^m$ is contact $CR$-product if and only if $P$ satisfies
$$ (\nabla_U P)V=g(PU_D,PV)\xi_\alpha-\eta^\alpha(V)U_D+\eta^\alpha(U_D)\eta^\alpha(V)\xi_\alpha, \leqno(2.21) $$
\noindent where $U,V$ tangent to $M^m$ and we are taking $U_D$ as
the component of $D$.
\end{thm}
\noindent {\bf Proof.} Since $\phi\equiv P$ on $N^T$, due to
indefinite $S$-structure of $\tilde{M}^{2n+s}$ using the Gauss
formula we get
$$(\nabla_XP)Y=g(PX,PY)\xi_\alpha-\eta^\alpha(Y)X+\eta^\alpha(X)\eta^\alpha(Y)\xi_\alpha+h(X,PY)-\phi h(X,Y), $$
\noindent for any $X,Y\in N^T $. Taking the components in $D$ one
gets
$$(\nabla_XP)Y=g(PX,PY)\xi_\alpha-\eta^\alpha(Y)X+\eta^\alpha(X)\eta^\alpha(Y)\xi_\alpha. \leqno(2.22) $$

\noindent Consider now $Z\in N^\perp$ and $Y\in N^T$. Similarly, we
can prove
$$(\nabla_Z P)Y=-\eta^\alpha(Y)Z, \leqno(2.23) $$
\noindent as consequence
$$h(Z,PY)=\phi h(Z,Y)+\eta(Y)Z,~~~~ \forall~ Y\in N^T,~~Z\in N^\perp. $$
\noindent Now it is easy to show that $(\nabla_U P)Z=0$ for all
$U\in \chi(M)$, $Z\in D^\perp$ and hence the conclusion.\\

\noindent Conversely, consider (2.21) exists. Let $U=X$, $V=Z$ with
$X\in D$ and $Z\in D^\perp$. The relation (2.21) becomes $(\nabla_X
P)Z=0$ and by using (2.6) we obtain $th(X,Z)=-A_{FZ} X$. Considering
$U=Z,V=X$ (with $X,Z$ as above) we obtain $(\nabla_Z
P)X=-\eta^\alpha (X)Z$. Thus one gets
$$A_{FZ}X=\eta^\alpha(X)Z, \leqno(2.24) $$
\noindent for all $X\in D$ and $Z\in D^\perp$. After the
computations we obtain $\tilde{g}(h(X,PY)-h(PX,Y))=0$. Let $X\in
H(M)$, $Z,W\in D^\perp$. Due to (2.24) we have
$$ \tilde{g}(h(X,Z),\phi W)=\tilde{g}(A_{FW} X,Z)=g(\eta^\alpha(X)W,Z)=\eta^\alpha (X)g(W,Z)=0. $$

\noindent Thus by virtue of Proposition (2.1), $N^\perp$ is totally
geodesic in $M^m$. Let now $X,Y\in D$, from (2.21) and (2.6) we
obtain $th(X,Y)=0$. If $Z\in D^\perp$ we have
$$0=\tilde{g}(th(X,Y),Z)=\tilde{g}(\tilde{\nabla}_X Y,\phi Z)=-\tilde{g}(Y,(\tilde{\nabla}_X \phi)Z)-\tilde{g}(Y,\phi\tilde{\nabla}_X Z)=-\tilde{g}(\phi Y,\nabla_X Z),  $$

\noindent replacing $Y$ by $\phi Y$ one obtains
$\tilde{g}(Y,\nabla_X Z)=0$ for all $X,Y \in D$ and $Z\in D^\perp$.
It follows that $g(\nabla_X Y,Z)=0$ which means that $N^T$ is also
totally geodesic in $M^m$. We may conclude that $M^m$ is a contact
$CR$-product in $\tilde{M}^{2n+s}$.\\

\begin{pp} Let $M^m$ be a contact $CR$-submanifold of an indefinite
$S$-manifold $\tilde{M}^{2n+s}$ with $\xi_\alpha \in D$. Then $M^m$
is a contact $CR$-product if and only if
$$A_{\phi Z}X=\eta^\alpha (X)Z, \leqno(2.25) $$
\noindent for all $X\in D$ and $Z\in D^\perp$.
\end{pp}
\noindent{\bf Proof.} Suppose that (2.25) holds. We have
$$\tilde{g}(B(X,Z),\phi W) = \tilde{g}(A_{\phi W}X,Z) = \eta^\alpha(X)g(Z,W) = 0,~~~~ \forall X \in H(M),~~ \forall Z,W \in D^\perp. $$

\noindent From Proposition 2.1 we get that $N^\perp$ (the integral
manifold of $D^\perp$) is totally geodesic in $M^m$. Consider now
$X,Y \in D$ and $Z \in D^\perp$. We have

$$ \tilde{g}(B(X,\phi Y ),\phi Z) = \tilde{g}(A_{\phi Z}X, \phi Y ) = \tilde{g}(\eta^\alpha(X)Z, \phi Y ) = 0. $$

\noindent Similarly $\tilde{g}(B(Y,\phi X),\phi Z) = 0$ and by Lemma
2.3 it follows that $D$ is completely integrable. To prove that
$N^T$ (the integral manifold of $D$) is totally geodesic in $M^m$ we
will prove that $\nabla_XY$ belongs to $N^T$ for all $X,Y$ tangent
to $N^T$. We have $g(\nabla_XY,Z) = -g(Y,\nabla_X Z)$. On the other
hand, from the hypothesis $\tilde{g}(B(X,Y),\phi Z)=0$. Then
$$ \tilde{g}(B(X,Y),\phi Z) = -\tilde{g}(Y,\tilde{\nabla}_X(\phi Z))=-\tilde{g}(\phi Y,\tilde{\nabla}_X Z)=-\tilde g(\phi Y,\nabla_X Z).$$

\noindent So, we obtain $g(\phi Y,\nabla_X Z)=0$, for all $X, Y \in
D$ and for all $Z\in D^\perp$. But $g(\xi_\alpha,\nabla_X Z)=0$ and
hence $g(Y,\nabla_XZ) =0$. We may conclude now that $\nabla_X Y \in
N^T$ for all $X,Y \in N^T$. Therefore the two integral manifolds
$N^T$ and $N^\perp$ are both totally geodesic in $M^m$ Consequently,
$M^m$ is locally a Riemannian product of $N^T$
and $N^\perp$.\\

\noindent Conversely, from the totally geodesy of $N^T$ and
$N^\perp$, using the Gauss formula we get
$\tilde{g}(\tilde{\nabla}_XY,\phi Z)=\tilde{g}(B(X,Y ), \phi Z)$
with $X,Y\in D$ and $Z\in D^\perp$. The right side is exactly
$g(A_{\phi Z}X, Y)$ while the left side equals to $-\tilde{g}(\phi
Y,\tilde{\nabla}_X Z)=g(\nabla_X(\phi
Y),Z)=0$.\\

\noindent It follows that $A_{\phi Z}X \in D^\perp$. Again by using
the Gauss formula we obtain after the computations
$\eta^\alpha(X)g(Z,W) = \tilde{g}(A_{\phi Z}X,W)$. Taking into
account that $A_{\phi Z}X \in D^\perp$ it follows $A_{\phi Z}X =
\eta^\alpha(X)Z$. This completes the proof.\\

\noindent Now we will prove the geometrical description of contact $CR$-products in indefinite $S$-space forms.\\

\begin{thm} Let $M^m$ be a generic, simply connected contact $CR$-submanifold of an indefinite $S$-space form
$\tilde{M}^{2n+s}(c)$, If $M^m$ is a contact $CR$-product then
\begin{enumerate}
 \item[{(i)}] For spacelike global vector field $\xi_\alpha$
$$(\frac{c+3}{4})(g(PX,PZ)FY-g(PY,PZ)FX)=0.$$
\item[{(ii)}] For timelike global vector field $\xi_\alpha$
 $$(\frac{c+5}{4})(g(PX,PZ)FY-g(PY,PZ)FX)=0. $$
\end{enumerate}

\noindent Therefore if $c\neq-3,-5$ and $M^m$ is a
$\phi$-anti-invariant submanifold of $\tilde{M}^{2n+s}$, then $M^m$
is locally a Riemannian product of an integral curve of $\xi_\alpha$
and a totally real submanifold $N^\perp$ of $\tilde{M}^{2n+s}$ and
if $c=-3,-5$, then $M^m$ is locally a Riemannian product of $N^T$
and $N^\perp$.
\end{thm}

\noindent {\bf Proof.} Since $M^m$ is generic it follows that
$\xi_\alpha \in D$. By remark 2.1 we have $h(X,Y)=0 $ for all $X,Y
\in D$ and $A_{FZ} X=\eta^\alpha (X) Z$ for all $X\in D$ and $Z\in
D^\perp$. Sice $T(M)^\perp=\phi D^\perp $ and $h\in T(M)^\perp$ by
using the Weingarten formula we immediately see that $g(h(X,Z),\phi
W)=g(A_{\phi W}X,Z)=\eta^\alpha(X)g(W,Z)$.\\

\noindent Consequently $h(X,Z)=\eta^\alpha(X)\phi Z$ for all $X\in D$ and $Z\in D^\perp$.\\

\noindent By making use of (2.11) we obtain for $X,U,V \in T(M)$.
$$(\nabla_X h)(U,PV)=-h(U,(\nabla_X P)V+P\nabla_X V)~~~~~~~~~~~~~~~~~~~~~~~~~~~~~~~~~~~~~~~~ \leqno(2.26)  $$
$$~~~~~~~~~~~~~~~~~~~~~=g(U,g(PX,PV)\xi_\alpha-\eta^\alpha(V) X+\eta^\alpha(X)\eta^\alpha(V)\xi_\alpha) $$
$$~~~~~~~~~~~~~~~~~~~~~~~~~~~~~~~~~~~~~~=g(PX,PV)\eta^\alpha(U)-\eta^\alpha(V)g(U,X)+\eta^\alpha(X)\eta^\alpha(V)\eta^\alpha(U), $$

\noindent hence we get
$$(\nabla_X h)(U,PV)=g(PX,PV)FU-\eta^\alpha(V)g(U,X)+\eta^\alpha(X)\eta^\alpha(V)\eta^\alpha(U) \leqno(2.27) $$

\noindent Substitute in (2.12) $Z$ by $PZ$ (with $Z\in T(M)$
arbitrary) the following identity holds:
$$(\nabla_X h)(Y,PZ)-(\nabla_Y h)(X,PZ)=\frac{c-\epsilon}{4} \{g(PX,PZ)FY-g(PY,PZ)FX \},  $$

\noindent where $\epsilon=\pm 1$ according to whether $\xi_\alpha$
is spacelike or timelike. Combining with (2.27) above relation
yields to
$$g(PX,PZ)FY-g(PY,PZ)FX =\frac{c-\epsilon}{4} \{g(PX,PZ)FY-g(PY,PZ)FX \},  $$

\noindent which is equivalent to
$$(\frac{c-\epsilon}{4}+1)(g(PX,PZ)FY-g(PY,PZ)FX)=0. \leqno(2.28) $$
\noindent Now we have to discuss two situations: $\epsilon=\pm 1$.
For spacelike global vector field $\epsilon=+1$, above equation
becomes
$$(\frac{c+3}{4})(g(PX,PZ)FY-g(PY,PZ)FX)=0. \leqno(2.29) $$

\noindent For timelike global vector field $\epsilon=-1$, the above
equation becomes
$$(\frac{c+5}{4})(g(PX,PZ)FY-g(PY,PZ)FX)=0. \leqno(2.30) $$
\noindent Now we discuss the two cases\\

\noindent {\bf Case I.} For $c\neq-3,-5,$ From the equation (2.29)
and (2.30) we obtain $g(PY,PZ)FX-g(PX, PZ)FY = 0$, for all $X, Y,Z
\in T(M)$: Since $M^m$ is generic we have $F \neq 0$ and it is not
difficult to prove that $P = 0$, thus $M^m$ is
$\phi$-anti-invariant. Moreover, by Theorem 2.6 we can say that
$M^m$ is a contact $CR$-product of an integral curve of $\xi_\alpha$
and a totally real submanifold $N^\perp$ of $\tilde{M}^{2n+s}$.\\

\noindent {\bf Case II.} For $c=-3,-5$ from (2.29) and (2.30). $M$
is a contact $CR$-product of the invariant submanifold $N^T$ and the
anti-invariant submanifold $N^\perp$: Since $N^T$ is totally
geodesic in $M^m$ and $h(X, Y)=0$ for all $X,Y\in D$ then $N^T$ is
totally geodesic in $\tilde{M}^{2n+s}$. Thus, we can use the well
known result that $M^m$ has constant $\phi-$ sectional curvature,
then $M^m$ is simply connected and hence $M^m$ is the Riemannian product of $N^T$ and $N^\perp$.\\

\noindent Let $\tilde{H}_h(U, V )$ be the $\phi$-holomorphic
bisectional curvature of the plane $U\wedge V$, i.e.
$$\tilde{H}_h(U, V )=\tilde{R}(\phi U,U.\phi V,V)~~~~~for~~U,V\in T(M). $$
\noindent We prove the following important lemmas for later use.\\

\begin{lm} Let $M^m$ be a contact $CR$-product of a indefinite $S$-manifold $\tilde{M}^{2n+s}$. Then, for any unit vector fields $X\in D$ and $Z\in
D^\perp$, then:
$$\tilde{g}(h(\nabla_{\phi X}X,\phi Z),Z)=-s,~~~\tilde{g}(h(\nabla_X {\phi X},\phi Z),Z)=s, \leqno(2.31) $$
$$\tilde{g}(h(X,\nabla_{\phi X}\phi Z),Z)=0,~~~\tilde{g}(h(\phi X,\nabla_X\phi Z),Z)=0. \leqno(2.32)$$
\end{lm}
\noindent {\bf Proof.} By theorem (2.6)
$$\tilde{g}(h(\nabla_{\phi X}X,\phi Z),Z)=\sum\limits_{\alpha=1}\limits^s \eta^\alpha(\nabla_{\phi X} X)g(\phi^2 Z,Z)=-\sum\limits_{\alpha=1}\limits^s g(X,\nabla_{\phi X} \xi_\alpha)g(Z,Z) $$
$$=-\sum\limits_{\alpha=1}\limits^s g(\phi X, \phi X)=-s,~~~~~~~~~~~~~~$$
\noindent which is first part of (2.31). We also have
$$\tilde{g}(h(\nabla_X {\phi X},\phi Z),Z)=\sum\limits_{\alpha=1}\limits^s \eta^\alpha(\nabla_X {\phi X} )g(Z,Z)=-\sum\limits_{\alpha=1}\limits^s g(\phi X,\nabla_X \xi_\alpha)g(Z,Z) $$
$$=\sum\limits_{\alpha=1}\limits^s g(\phi X, \phi X)=s,~~~~~~~~~~~~~~~~$$
\noindent which is other part of (2.31). Finally
$$\tilde{g}(h(X,\nabla_{\phi X}\phi Z),Z)=\sum\limits_{\alpha=1}\limits^s \eta^\alpha(X)g(\nabla_{\phi X} \phi Z,Z)=0,$$
$$\tilde{g}(h(\phi X,\nabla_{X}\phi Z),Z)=\sum\limits_{\alpha=1}\limits^s \eta^\alpha(\phi X)g(\nabla_{X} \phi Z,Z)=0,$$
\noindent which are (2.32).

\begin{lm} Let $M^m$ be a contact $CR$-product of a indefinite $S$-manifold $\tilde{M}^{2n+s}$. Then, for any unit vector fields $X\in D$ and $Z\in D^\perp$ we have
$$\tilde{H}_h(X,Z)=2s-2\|h(X,Z)\|^2  \leqno(2.33) $$
\end{lm}
\noindent {\bf Proof.} We know
$$\tilde{R}(\phi X,X,\phi Z,Z)=\tilde{g}(({\nabla^h}_{\phi X}h)(X,\phi Z)-({\nabla^h}_{X}h)(\phi X,\phi Z),Z)  $$
\noindent by using (2.11) and above lemma, we gets
$$\tilde{R}(\phi X,X,\phi Z,Z)=\tilde{g}(\tilde{\nabla}_{\phi X} h(X,\phi Z)-h(\nabla_{\phi X}X,\phi Z)-h(X,\nabla_{\phi X}\phi Z),Z)$$
$$~~~~~~~~~~~~~~~~~~~~~~-\tilde{g}(\tilde{\nabla}_{X} h(\phi X,\phi Z)-h(\nabla_{X}\phi X,\phi Z)-h(\phi X,\nabla_{X}\phi Z),Z) $$
$$~~~~~~~~~~~~~~~~~~~~=2s-\tilde{g}(h(X,\phi Z),\nabla_{\phi X}Z)-\tilde{g}(h(X,\phi Z),h(\phi X,Z)) $$
$$~~~~~~~~~~~~~~~~~~~~+g(h(\phi X,\phi Z),\nabla_X Z)+g(h(\phi X,\phi Z),h(X,Z)).  $$
\noindent Now, using the fact that $\nabla_X Z$ and $\nabla_{\phi X}
Z$ belong to $D^\perp$, therefore above equation becomes
$$\tilde{R}(\phi X,X,\phi Z,Z)=2s-2\|h(X,Z\|^2, $$
\noindent this ends the proof.\\

\noindent We come to know that $\tilde{H}_h(U,\xi_\alpha)=0$ and
$h(U,\xi_\alpha)=\phi U$, So, when we will refer to the
$\phi$-holomorphic bisectional curvature of the plane $U\wedge V$,
we intend that this plane is orthogonal $\xi_\alpha$. Thus for $X$
in the above lemma we can suppose that it belongs to $H(M)$.\\

\begin{pp} Let $\tilde{M}^{2n+s}(c)$ be a indefinite $S$-space form and let $X,Z$ be two unit vector fields
orthogonal to global vector filed $\xi_\alpha$. Then the
$\phi$-holomorphic bisectional curvature of the plane $X\wedge Z$ is
given by
\begin{enumerate}
 \item[{(i)}] For spacelike global vector field $\xi_\alpha$
 $$ \tilde{H}_B (X,Z)=\frac{c-1}{4}g(\phi X,Z)-\frac{c+3}{4}g(\phi X,Z)^2+c. $$
 \item[{(ii)}] For timelike global vector field $\xi_\alpha$
 $$ \tilde{H}_B (X,Z)=\frac{c+1}{4}g(\phi X,Z)-\frac{c-3}{4}g(\phi X,Z)^2+c. $$
\end{enumerate}
\end{pp}

\begin{cor}
Let $\tilde{M}^{2n+s}(c)$ be a indefinite $S$-space form and let
$X\in H(M)$ and $Z\in D^\perp$ be unit vector fields orthogonal to
global vector filed $\xi_\alpha$. Then the $\phi$-holomorphic
bisectional curvature of the plane $X\wedge Z$ for spacelike and
timelike global vector field $\xi_\alpha$ is given by
$$ \tilde{H}_B (X,Z)=c.$$
\end{cor}

\begin{thm} Let $\tilde{M}^{2m+s}(c)$ be a indefinite $S$-space form
and let $ M=N^T\times N^\perp$ be a contact $CR$-product in
$\tilde{M}^{2n+s}$. Then the norm of the second fundamental form of
$M$ satisfies the inequality
$$ \|h\|^2\geq((3c+8s-3\epsilon)p+2s)q,  $$
\noindent where $\epsilon=\pm1$ according to whether global vector
field $\xi_\alpha$ is spacelike or timelike. The equality sign holds
if and only if both $N^T$ and $N^\perp$ are totally geodesic.
\end{thm}
\noindent {\bf{Proof.}} For $X\in H(M)$ and $Z\in D^\perp$ we have
$\|h(X,Z)\|^2=\frac{1}{4}(3c+8s-3\epsilon)$.\\

\noindent Now, we choose a local field of orthonormal frames
$$ \{X_1,...,X_p,X_{p+1}=\phi X_1,...,X_{2p}=\phi X_p,X_{2p+1}=Z_1,...,X_n=Z_q,~~~~$$
$$~~~~~~~~X_{n+1}=\phi Z_1,...,X_{n+q}=\phi Z_q, X_{n+q+1},...,X_2m,\xi_1,...,\xi_s \} $$
\noindent on $\tilde{M}^{2n+s}(c)$ in such a way that
$\{X_1,...,X_{2p} \}$ is a local frame field on $D$ and
$\{Z_1,...,Z_q \}$ is a local frame field on $D^\perp$. Thus
$$ \|h\|^2=\|h(D,D)\|^2+\|h(D^\perp,D^\perp\|^2+2\|h(D,D^\perp\|^2\geq 2\|h(D,D^\perp\|^2 $$
$$ =2(\sum\limits_{i=1}\limits^{2p}\sum\limits_{j=1}\limits^q\|h(X_i,Z_j)\|^2+\sum\limits_{\alpha=1}\limits^s\sum\limits_{j=1}\limits^q \|h(\xi_\alpha,Z_j\|^2) $$
$$=((3c+8s-3\epsilon)p+2s)q, $$
\noindent where $X_i$ and $Z_j$ are orthonormal basis in $H(M)$ and
$D^\perp$ respectively. The equality sign holds if and only if
$h(D,D)=0$ and $h(D^\perp,D^\perp)= 0$, which is equivalent to the
totally geodesy of $N^T$ and $N^\perp$.

\section{\large $CR$-warped product submanifolds in indefinite $S$-manifolds}

\noindent The main purpose of this section is devoted to the
presentation of some properties of warped product contact $CR$-submanifolds in indefinite $S$-manifolds.\\

\noindent Let $(N_1, g_1)$ and $(N_2, g_2)$ be two Riemannian
manifolds and $f$, a positive differentiable function on $N_1$. The
{\it{warped product}} of $N_1$ and $N_2$ is the product manifold
$N_1\times{_{f}N_2}=(N_1\times N_2,~g)$, where
$$g=g_1+f^2g_2,~~~\leqno(3.1)$$
where $f$ is called the {\it{warping function}} of the warped product.
The warped product $N_1\times{_{f} N_2}$ is said to be {\it{trivial}} or
simply Riemannian product if the warping function $f$ is constant. This means that the Riemannian product is a special case of warped product.\\

\noindent
We recall the following general results obtained by Bishop and O'Neill [2] for warped product manifolds.\\

\begin{lm} Let $M=N_1\times{_{f} N_2}$ be a warped product manifold with
the warping function $f$, then
\begin{enumerate}
\item [{(i)}] $\nabla_XY\in TN_1$ {\it{for each}} $X, Y\in TN_1$,
\item [{(ii)}] $\nabla_XZ=\nabla_ZX=(X\ln f)Z$, {\it{for each $X\in TN_1$ and}} $Z\in TN_2$,
\item [{(iii)}] $\nabla_ZW=\nabla_Z^{N_2}W-\frac{g(Z,W)}{f}{\mbox{{\it{grad}}}}f$,
\end{enumerate}
\noindent where $\nabla$ and $\nabla^{N_2}$ denote the Levi-Civita
connections on $M$ and $N_2$, respectively.
\end{lm}

\noindent In the above lemma {\it{grad}}$f$ is the gradient of the
function $f$ defined by $g({\mbox{{\it{grad}}}}f, U)=Uf$, for each
$U\in TM$. From the Lemma 3.1, we have on a warped product manifold
$M=N_1\times{_{f} N_2}$
\begin{enumerate}
\item [{(i)}] $N_1$ is totally geodesic in $M$.
\item [{(ii)}] $N_2$ is totally umbilical in $M$.
\end{enumerate}

\begin{thm} Let $\tilde{M}^{2n+s}$ be a indefinite $S$-manifold and let $M = N^\perp \times_f N^T$ be a warped product $CR$-submanifold
such that $N^\perp$ is a totally real submanifold and $N^T$ is
$\phi$-holomorphic (invariant) of $\tilde{M}^{2n+s}$. Then $M$ is a
$CR$-product.
\end{thm}

\noindent {\bf Proof.} Let $X$ be tangent to $N^T$ and let $Z$ be a
vector field tangent to $N^\perp$. From the above lemma we
find that
$$\nabla_XZ = (Z ln f) X. \leqno(3.2)) $$

\noindent Now we have two cases either $\xi_\alpha$ is tangent to $N^T$ or $\xi_\alpha$ is tangent to $N^\perp$.\\

\noindent Case I. $\xi_\alpha$ is tangent to $N^T$. Take $X
=\xi_\alpha$. Since $\nabla_Z \xi_\alpha =-PZ = 0$ and $\nabla_Z
{\xi_\alpha} = \nabla_{\xi_\alpha} Z$ ($\xi_\alpha$ is tangent to
$N^T$ while $Z$ is tangent to $N^\perp$) one gets $0=Z(ln
f)\xi_\alpha$ and hence $Z(ln f)=0$ for all $Z$ tangent to
$N^\perp$. Consequently $f$ is constant and thus the warped product
above is nothing but a Riemannian product. \\

\noindent Case II. Now we will consider the other case, $\xi_\alpha$
is tangent to $N^\perp$. Similarly, take $Z =\xi_\alpha$. Since
$\nabla_X \xi_\alpha =-\epsilon_\alpha PX =
-\epsilon_\alpha \phi X$ it follows $-\epsilon_\alpha\phi X = (\xi_\alpha ln f)X$. But this is impossible if $dim N^T \neq 0$.\\

\noindent {\bf Remark:}  There do not exist warped product
$CR$-submanifolds in the form $N^\perp \times_f N^T$ other than
$CR$-products such that $N^T$ is a $\phi$-invariant submanifold and
$N^\perp$ is a totally real submanifold of $\tilde{M}$.\\

\noindent From now on we will consider warped product
$CR$-submanifolds in the form $N^T \times_f N^\perp$.\\

We can say that {\it A contact $CR$-submanifold $M$ of a indefinite
$S$-manifold $\tilde{M}^{2n+s}$, tangent to the structure vector
field $\xi_\alpha$ is called a contact $CR$-warped product if it is
the warped product $N^T \times_f N^\perp$ of an invariant
submanifold $N^T$, tangent to $\xi_\alpha$ and a totally real
submanifold $N^\perp$ of $\tilde{M}^{2n+s}$, where $f$ is the warping function}.\\

\begin{lm} Let $M^m$ be a contact $CR$-submanifold in indefinite $S$-manifold $\tilde{M}^{2n+s}$, such that $\xi_\alpha \in D$. Then we have
$$ {g}(\nabla_U X,Z)=-\tilde{g}(\phi A_{\phi Z} U,X),~~~~\forall~ X\in D,~~\forall~ Z\in D^\perp,~~\forall~ U\in T(M), \leqno(3.3)  $$
$$ A_{\phi \mu}X+A_{\mu}\phi X=0~~~~\forall~ X\in D,~~\forall~ \mu\in\nu.~~~~~~~~~~~~~~~~~~~~~~~~~~~~ \leqno(3.4) $$
\end{lm}
\noindent {\bf Proof.} We have $\tilde{g}(\phi A_{\phi Z}
U,X)=\tilde{g}(A_{\phi Z}U,\phi X)=\tilde{g}({\nabla^\perp}_U \phi
Z-\tilde{\nabla}_U \phi Z,\phi X)=-\tilde{g}(\phi \tilde{\nabla}_U
Z,\phi X)=-\tilde{g}(\tilde{\nabla}_U
Z,X)+\eta^\alpha(\tilde{\nabla}_UZ)\eta^\alpha(X)=-g(\nabla_UZ,X)-\tilde{g}(Z,\tilde{\nabla}_U
\xi)\eta^\alpha(X)=-g(\nabla_U Z,X)+g(Z,\phi
U)\eta^\alpha(X)=-g(\nabla_U Z,X) $. So, equation (3.3) follows.\\

\noindent For the proof of the equation (3.4) we have $g(A_{\phi
\mu}X,U)=-g(\tilde{\nabla}_X
\phi\mu,U)=g(\mu,\phi\tilde{\nabla}_XU)$ and $g(A_{\mu} U,\phi
X)=-g(\mu,\phi\tilde{\nabla}_U X)$ with $U\in T(M)$. It follow that
$A_{\phi \mu}X+A_{\mu}\phi X=0$, $\forall~ X\in D$, $\forall~\mu\in\nu$.\\

\begin{lm} If $M=N^T\times_f N^\perp$ is a contact $CR$-warped product
in a indefinite $S$-manifold $\tilde{M}^{2n+s}$ the for $X$ tangent
to $N^T$ and $Z,W$ tangent to $N^\perp$ we have
$$ g(h(D,D^\perp),\phi D^\perp)=0 \leqno(3.5) $$
$$ \xi_\alpha(f)=0 \leqno (3.6)$$
$$ g(h(\phi X,Z),\phi W)=(X ln f)g(Z,W) \leqno(3.7) $$
\end{lm}
\noindent {\bf Proof.} For any $X,Y\in D$ and $Z\in D^\perp$, we
have
$$g(h(X,Y),\phi Z)=g(\tilde{\nabla}_X Y,\phi Z)=-g(\phi Y,\tilde{\nabla}_X Z)=g(\tilde{\nabla}_X \phi Y,Z)=0.  $$

\noindent We know that $\nabla_U \xi_\alpha=-\epsilon_\alpha PU$. It
follows that $\nabla_Z \xi_\alpha=0$ for all $Z$ tangent to
$N^\perp$. Using lemma 3.1 and Theorem 3.2 we get equation (3.6).\\

\noindent From equation (3.2) it follows that
$$g(h(\phi X,Z),\phi W)=g(A_{\phi W}Z,\phi X)=-g(\nabla_Z W,X)=X(ln f) g(Z,W). $$
\noindent Hence proved.\\

\begin{thm} The necessary and sufficient condition for a strictly proper $CR$-submanifold
$M$ of a indefinite $S$-manifold $\tilde{M}^{2n+s}$, tangent to the
structure vector field $\xi_\alpha$ to be locally a contact
$CR$-warped product is that
$$ A_{\phi Z} X=(-(\phi X)(\mu)-\eta^\alpha(X))Z+\eta^\alpha(X)\eta^\alpha(Z)\xi_\alpha \leqno(3.7) $$
\noindent for some function $\mu$ on $M$ satisfying $W\mu=0$ for all
$W\in D^\perp$.
\end{thm}
\noindent {\bf Proof.} Let $M=N^T\times_fN^\perp$ be a locally
contact $CR$-warped product.\\

\noindent Consider $X,Y\in D$, $Z\in D^\perp$. We can easily get
$g(A_{\phi Z}X,Y)=0$, which shows that $A_{\phi Z}X$ belongs to
$D^\perp$.\\

\noindent Now take any $W\in D^\perp$, we get
$$ g(A_{\phi Z}X,W)=(-(\phi X)(\mu)-\eta^\alpha(X))g(W,Z)+\eta^\alpha(X)\eta^\alpha(Z)\eta^\alpha(W) $$
\noindent hence the result where $\mu=ln f$.\\

\noindent Conversely, Let
$$A_{\phi Z} X=(-(\phi X)(\mu)-\eta^\alpha(X))Z+\eta^\alpha(X)\eta^\alpha(Z)\xi_\alpha. $$

\noindent We get easily that $\tilde{g}(h(\phi X,Y),\phi Z)=0$,\\

\noindent Also
$$\tilde{g}(h(X,W),\phi Z)=(-(\phi X)(\mu)-\eta^\alpha(X))Z+\eta^\alpha(X)\eta^\alpha(Z)\xi_\alpha,$$
\noindent where $X,Y\in D$ and $Z,W\in D^\perp$. In the above
equation replacing $X$ by $\phi X$, we obtain
$$\tilde{g}(h(\phi X,W),\phi Z)=(X(\mu)-\eta^\alpha(X)\xi_\alpha(\mu))g(Z,W). $$
\noindent So, if we take $X\in H(M)$ it becomes $\tilde{g}(h(\phi
X,W),\phi Z)=(X(\mu)g(Z,W)$ and if $X=\xi_\alpha$ we get
$\tilde{g}(h(\phi X,W),\phi Z)=(1-{\delta_\alpha}^\alpha)\xi_\alpha(\mu)g(Z,W)$.\\

\section{\large A geometric inequality for contact $CR$-warped product in indefinite $S$-space form}

\noindent In this section, we will prove the main theorems of the
paper.\\

\noindent Let $M$ be a (pseudo-)Riemannian $k$-manifold with inner
product $<,>$ and $e_1,...,e_k$ be an orthonormal frame fields on
$M$. For differentiable function $\phi$ on $M$, the gradient $\nabla
\phi$ and the Laplacian $\triangle \phi$ of $\phi$ are defined
respectively by
$$ <\nabla \phi,X>=X(\phi), \leqno(4.1)$$
$$ \triangle \phi=\sum\limits_{j=1}\limits^k\{(\nabla_{e_j} e_j)\phi-e_j e_j(\phi) \}=-div \nabla \phi \leqno(4.2) $$

\noindent for vector field $X$ tangent to $M$, where $\nabla$ is the
Riemannian connection on $M$. As consequence, we have
$$\|\nabla \phi\|^2=\sum\limits_{j=1}\limits^k (e_j(\phi))^2. \leqno(4.3) $$

\begin{thm} Let $M=N^T\times_f N^\perp$ be a contact $CR$-warped
product of a indefinite $S$-space form $\tilde{M}^{2n+s} (c)$. Then
the second fundamental form of $M$ satisfies the following
inequality

$$ \|h\|^2\geq p\{3\|\nabla ln f\|^2-\triangle ln f+(c+2)k+1 \}. \leqno(4.4)  $$

\end{thm}
\noindent {\bf Proof.} We have
$$ \|h(D,D^\perp)\|^2=\sum\limits_{j=1}\limits^{k}\sum\limits_{i=1}\limits^p \|h(X_j,Z_i)\|^2, \leqno(4.6) $$
\noindent where $X_j$ for $\{j=1,...,k\}$ and $Z_\alpha$ for
$\alpha=\{1,...,p\}$ are orthonormal frames on $N^T$ and $N^\perp$,
respectively. On $N^T$ we will consider a $\phi$-adapted orthonormal
frame, namely $\{e_j,\phi e_j,\xi_\alpha\}$, where $\{j=1,...,k\},\{\alpha=1,...,s\}$.\\

\noindent We have to evaluate $\|h(X,Z)\|^2$ with $X\in D$ and $Z\in
D^\perp$. The second fundamental form $h(X,Z)$ is normal to $M$ so,
it splits into two orthogonal components
$$h(X,Z)=h_{\phi D^\perp}(X,Z)+h_\nu (X,Z), \leqno(4.7) $$
\noindent where $h_{\phi D^\perp}(X,Z) \in \phi D^\perp$ and $h_\nu
(X,Z) \in \nu$. So
$$ \|h(X,Z)\|^2=\|h_{\phi D^\perp}(X,Z)\|^2+\|h_\nu (X,Z)\|^2. \leqno(4.8)  $$
\noindent If $X=\xi_\alpha$, we have $h(\xi_\alpha,Z)=-\phi Z$.
Hence
$$ h_{\phi D^\perp}(\xi_\alpha,Z)=-\phi Z,~~~~h_\nu (\xi_\alpha,Z)=0. \leqno(4.9) $$
\noindent Consider now $X\in H(M)$ and let's compute the norm of the
$\phi D^\perp$-component of $h(X,Z)$. We have
$$ \|h_{\phi D^\perp}(X,Z)\|^2=<h_{\phi D^\perp}(X,Z),h(X,Z)>. $$
\noindent By using relation (3.7), after the computations, we obtain
$$ \|h_{\phi D^\perp}(X,Z)\|^2=(\phi X(ln f))^2 \|Z\|^2. $$
\noindent So
$$ \|h_{\phi D^\perp}(e_j,Z_i)\|^2=(\phi e_j(ln f))^2,~~~~\|h_{\phi D^\perp}(\phi e_j,Z_i)\|^2=(e_j (ln f))^2. \leqno(4.10) $$
\noindent On the other hand, from (4.2) we have
$$ \|\nabla ln f\|^2=\sum\limits_{j=1}\limits^k (e_j(ln f))^2+\sum\limits_{j=1}\limits^k (\phi e_j(ln f))^2. \leqno(4.11) $$
\noindent Since $\xi_\alpha(ln f)=0$. Finally we can compute the
norm $\| h_{\phi D^\perp}(D,D^\perp)\|^2$. Thus $\|h_{\phi
D^\perp}(D,D^\perp)\|^2=\sum\limits_{j=1}\limits^{k}\sum\limits_{i=1}\limits^p\{
\|h_{\phi D^\perp}(e_j, Z_i)\|^2+\|h_{\phi D^\perp}(\phi
e_j,Z_i)\|^2\}+\sum\limits_{\alpha=1}\limits^s\sum\limits_{i=1}\limits^p
\|h_{\phi D^\perp}(\xi_\alpha,Z_i)\|^2=\sum\limits_{i=1}\limits^p
\|\nabla ln f\|^2+\sum\limits_{i=1}\limits^p \|\phi Z_i\|^2. $
\noindent Since $\|\phi Z_i\|^2=1$ we can conclude that

$$ \|h_{\phi D^\perp}(D,D^\perp)\|^2=\{\sum\limits_{i=1}\limits^p \|\nabla ln f\|^2+1\}p. \leqno(4.12)  $$

\noindent Now we will compute the norm of the $\nu$-component of
$h(X,Z)$. We have

$$ \|h_\nu(X,Z)\|^2=<h_\nu(X,Z),h(X,Z)>=<A_{h_\nu (X,Z)} X,Z>, $$

\noindent by using lemma (3.3) we can write $A_{h_\nu(X,Z)}X=A_{\phi
h_\nu(X,Z)}(\phi X) $ so,

$$ \|h_\nu(X,Z)\|^2=<\phi h(X,Z)-\phi h_{\phi D^\perp}(X,Z),h(\phi X,Z)>.  $$

\noindent Since $\phi h_{\phi D^\perp}(X,Z)$ belongs to $D^\perp$ we
obtain

$$ \|h_\nu(X,Z)\|^2=\tilde{g}(\phi h(X,Z),h(\phi X,Z)),~~~~\forall~ X\in H(M),~~Z\in D^\perp. \leqno(4.13)  $$

\noindent Consider the tensor field $\tilde{H}_B$. As we already
have seen
$$ \tilde{H}_B (X,Z)=<(\nabla_{\phi X} )h(X,Z)-(\nabla_X h)(\phi X,Z),\phi Z>,~~~~\forall~ X\in H(M),~~Z\in D^\perp. $$

\noindent Using the definition of $\nabla h$, developing the
expression above we obtain
$$ \tilde{H}_B (X,Z)=<{\nabla^\perp}_{\phi X} h(X,Z),\phi Z>-<h(\nabla_{\phi X} X,Z),\phi Z>-<h(X,\nabla_{\phi X} Z),\phi Z>$$
$$~~~~~~~~~~~~~~~~~~~-<{\nabla^\perp}_{X} h(\phi X,Z),\phi Z>+<h(\nabla_{X} \phi X,Z),\phi Z>+<h(\phi X,\nabla_{X} Z),\phi Z>, $$
\noindent after using lemma 3.4 and theorem 3.5, we get
$$ \tilde{H}_B (X,Z)=\|Z\|^2\{(\phi X(ln f))^2-(\phi X)(\phi X(ln f))+(X ln f)^2-X(X ln f) $$
$$~~~~~~~~~+ (\phi \nabla_{\phi X} X)(ln f)-\|X\|^2-(\phi X \nabla_X \phi X)(ln f)-\|X\|^2 $$
$$~~~~~~~~~~~+(\phi X (ln f))^2+(X (ln f))^2 \}+2<\phi h(X,Z),h(\phi X,Z)>, $$
\noindent which becomes
$$ \tilde{H}_B (X,Z)=\|Z\|^2\{2(\phi X(ln f))^2-(\phi X)^2(ln f)+2( X lnf)^2-X^2 (ln f) $$
$$~~~~~~~~~+ ((\phi \nabla_{\phi X} X)-(\phi X \nabla_X \phi X))(ln f)-2\|X\|^2\}+2\|h_{\nu} (X,Z)\|^2. \leqno(4.14)  $$
\noindent We can easily prove that
$$ (\phi \nabla_{\phi X}X)(ln f)=(\nabla_{\phi X}\phi X)(ln f),~~~~(\phi \nabla_X{\phi X})(ln f)=-\nabla_X X(ln f). \leqno(4.15)$$

\noindent Using (4.14) and (4.14), we get

$$ \tilde{H}_B (X,Z)=\|Z\|^2\{2(\phi X(ln f))^2-(\phi X)^2(ln f)+2( X lnf)^2-X^2 (ln f) \leqno(4.16) $$
$$~~~~~~~+ (\nabla_{\phi X}(\phi X)+\nabla_X X)(ln f)-2\|X\|^2\}+2\|h_{\nu} (X,Z)\|^2.   $$

\noindent Using orthonormal frames, we have

$$ \tilde{H}_B (e_j,Z_i)=\|Z_i\|^2\{2(\phi e_j(ln f))^2-(\phi e_j)^2(ln f)+2( e_j lnf)^2-{e_j}^2 (ln f) \leqno(4.17) $$
$$~~~~~~~~~+ (\nabla_{\phi e_j}(\phi e_j)+\nabla_{e_j} e_j)(ln f)-2\|X_j\|^2\}+2\|h_{\nu} (e_j,Z_i)\|^2.   $$
\noindent Similarly,
$$ \tilde{H}_B (\phi e_j,Z_i)=\|Z_i\|^2\{2( e_j(ln f))^2-{e_j}^2(ln f)+2(\phi e_j lnf)^2-(\phi e_j)^2 (ln f) \leqno(4.18) $$
$$~~~~~~~~~+ (\nabla_{e_j}e_j+\nabla_{\phi e_j} (\phi e_j))(ln f)-2\|X_j\|^2\}+2\|h_{\nu} (\phi e_j,Z_i)\|^2.   $$
\noindent On the other hand we have
$$ \triangle (ln f)=\sum\limits_{j=1}\limits^k\{(\nabla_{e_j} e_j)(ln f)-{e_j}^2(ln f) \}+\sum\limits_{j=1}\limits^k\{(\nabla_{\phi e_j} \phi e_j)(ln f)-{\phi e_j}^2(ln f) \} $$
\noindent by using (4.3), we get
$$ 2\|\nabla ln f\|^2=2\sum\limits_{j=1}\limits^k (e_j(ln f))^2+2\sum\limits_{j=1}\limits^k (\phi e_j(ln f))^2. $$

\noindent Since $\xi_\alpha ln f=0$. Taking the sum of (4.17) and
(4.18) we get
$$ 2\sum\limits_{j=1}\limits^k\sum\limits_{i=1}\limits^p\{ \|h_{\nu} (e_j,Z_i)\|^2+\|h_{\nu} (\phi e_j,Z_i)\|^2\}=\sum\limits_{j=1}\limits^k\sum\limits_{i=1}\limits^p\{\tilde{H}_B(e_j,Z_i)+\tilde{H}_B (\phi e_j,Z_i)\} $$
$$~~~~~~~~~~~~~~~~~~~~~~~~~~~~~~~~~~~~~~~~~~~~~~~~-2p\triangle (ln f)+4kp+4\|\nabla ln f\|^2, $$

\noindent by using proposition 2.11, we have
$$\sum\limits_{j=1}\limits^k\sum\limits_{i=1}\limits^p\{ \|h_{\nu} (e_j,Z_i)\|^2+\|h_{\nu} (\phi e_j,Z_i)\|^2\}=ckp-p\triangle (lnf)+2kp+2k\|\nabla ln f\|^2. \leqno(4.19) $$
\noindent Now from (4.8), (4.12) and (4.19) we conclude that $h$
satisfies the inequality.


\begin{thebibliography}{20}
\bibitem {} J. K. Beem, P. E. Ehrlich and K. Easley,{\it Global Lorentzian geometry}, Marcel Dekker, New York, 1996.
\bibitem {} J. K. Beem and P. E. Ehrlich,{\it Singularities, incompleteness, and the Lorentzian distance function}, Math. Proc. Camb. Phil. Soc. 85, 161 (1979).
\bibitem {} J. K. Beem, P. E. Ehrlich and T. G. Powell, {\it Warped product manifolds in relativity}, Selected studies : Physics-astrophysics, mathematics, history of science, pp.41–56, North-Holland, Amsterdam-New York, 1982.
\bibitem {} A. Bejancu, M. Kon and K. Yano,{\it CR submanifolds of a complex space form}, J. Diff. Geom., 16 (1981), 137-145.
\bibitem {} R. L. Bishop and B. O'Neill, {\it{Manifolds of negative curvature}}, Trans. Amer. Math. Soc. 145 (1969), 1-49.
\bibitem {} D. E. Blair and K. Yano, {\it{Affine almost contact manifolds and f-manifolds with affine Killing structure tensors}}, Kodai Math. Sem. Rep. 23 (1971) 473-479.
\bibitem {} D. E. Blair, {\it Riemannian geometry of contact and symplectic manifolds}, Progr. Math., 203, Birkh$\ddot{a}$user Boston, Boston, MA, 2002.
\bibitem {} V. Bonanzinga and K. Matsumoto, {\it{Warped product CR-submanifolds in locally conformal Kaehler manifolds}}, Periodica Math. Hungar. 48 (2004), 207-221.
\bibitem {} L. Brunetti and A. M. Pastore, {\it Curvature of a class of indefinite globally framed $f$-manifolds}, Bull. Math. Soc. Sci. Math. Roumanie (N.S.) 51(99)(2008), no. 3, 183-204.
\bibitem {} B. Y. Chen, {\it{Geometry of warped product CR-submanifolds in  Kaehler manifold}}, Monatsh. Math. 133 (2001), 177-195.
\bibitem {} B. Y. Chen, {\it{Geometry of warped product CR-submanifolds in Kaehler manifolds II}}, Monatsh. Math. 134 (2001), 103-119.
\bibitem {} S. I. Goldberg and K. Yano,{\it On normal globally framed $f-$manifolds}, T$\hat{o}$hoku Math. J. (22) (1970), 362-370.
\bibitem {} I. Hasegawa and I. Mihai, {\it Contact CR-warped product submanifolds in Sasakian manifolds}, Geom. Dedicata 102 (2003), 143-150.
\end{thebibliography}
\end{document}